\documentclass{amsart}
\usepackage{amsmath,amssymb,amsthm}

\def \({\left(}
\def \){\right)}
\def \M{\mathcal{M}}
\def \N{\mathcal{N}}
\def \H{\mathcal{H}}
\def \K{\mathcal{K}}
\def \S{\mathcal{S}}
\def \C{\mathbb{K}}
\def \A{\mathbb{A}}

\def \g{\mathfrak{g}}
\def \n{\mathbf{n}}

\newtheorem{lemma}{Lemma}

\newtheorem{thm}[lemma]{Theorem}
\theoremstyle{definition}

\theoremstyle{remark}
\newtheorem*{remark}{Remark}

\begin{document}
\title{On the variety of almost commuting nilpotent matrices}
\author{Eliana Zoque}
\address{Department of Mathematics, The University of Chicago}

\begin{abstract}
We study the variety of $n\times  n$ matrices with commutator of rank at most one. We describe its irreducible components; two of them correspond to the pairs of commuting matrices, and $n-2$ components of smaller dimension corresponding to the pairs of rank one commutator. In our proof we define a map to the zero fiber of the Hilbert scheme of points and study the image and the fibers.

\end{abstract}
\maketitle

\section{Introduction}

Let $V$ be a vector space of dimension $n$ over a field $\C$ of
characteristic equal to 0 or $\geq n/2$. Let
$\g=\mathfrak{gl}_n(V)$ and $\n$ be the nilcone of $\g$, i.e., the
cone of nilpotent matrices of $\g$. We write elements of $V$ and
$V^*$ as column and row vectors, respectively. In this paper we
study the variety
$$\N:=\{(X,Y,i,j)\in \n\times\n\times V\times
V^*\,|\,[X,Y]+ij=0\}$$ and prove that it has $n$ irreducible
components: 2 of dimension $n^2+n-1$ corresponding to the case
where the matrices commute, and $n-2$ of dimension $n^2+n-2$,
corresponding to the noncommutative pairs.

The pairs of almost commuting matrices have been studied recently
in \cite{Gan-Ginzburg:06}, where Gan and Ginzburg study the
structure of the scheme
$$\M:=\{(X,Y,i,j)\in \g\times \g\times V\times
V^*\,|\,[X,Y]+ij=0\}.$$ They prove that the irreducible components
of $\M$ are the closures of the sets $\M_0,\,\M_1,\,\dots,\,\M_n$,
defined as \begin{multline*}\M_t=\{(X,Y,i,j)\in\M\,|\,Y \text{ has
pairwise distinct eigenvalues and }\\\dim\C\langle X,Y\rangle
i=t,\, \dim j\C\langle X,Y\rangle=n-t\}\end{multline*} where
$\C\langle X,Y\rangle i$ (resp.\ $j\C\langle X,Y\rangle$) is the
smallest subspace of $V$ (resp.\ $V^*$) containing $i$ (resp.\
$j$) and invariant under $X$ and $Y$.

Let $\mathcal{K}_n=\{(X,Y)\in \n\times\n\,|\,[X,Y]=0\},$ the
variety of commuting nilpotent matrices. Baranovsky proved in
\cite{Baranovsky:01} that $\mathcal{K}_n$ is irreducible and has
dimension $n^2-1$. In his proof, he shows that
$U=\{(X,Y,i)\in\mathcal{K}_n\times V\,|\,\C[X,Y]i=V\}$ is
irreducible, dense in $\mathcal{K}_n\times V$ an has dimension
$n^2+n-1$. $GL(V)$ acts faithfully on $U$, and the quotient
$U/GL(V)$ is a fiber of the Hilbert scheme of points under the
Hilbert-Chow morphism.

Let $$\N_{r,s}=\{(X,Y,i,j)\in\N\,|\,\dim \C\langle X,Y\rangle
i=r,\,\dim j\C\langle X,Y\rangle=s \},$$
$$\N'_{r,s}=\{(X,Y,i,j)\in\N\,|\,\dim \C\langle X,Y\rangle i\leq
r,\,\dim j\C\langle X,Y\rangle\leq s \}$$ and $\overline\N_{r,s}$
the Zariski closure of $\N_{r,s}$. Clearly
$\N_{r,s}\subseteq\overline\N_{r,s}\subseteq\N'_{r,s}.$

Since $X$ and $Y$ can be put in upper triangular form
simultaneously (\cite{Etingof-Ginzburg:02}, Lemma 12.7), we have
that $[X,Y]$ is not only strictly upper triangular, but the
entries that are located two positions above the diagonal are zero
too. Therefore,
$$\N=\N_{0,n}\cup\N_{n,0}\cup\bigcup_{0<r+s<n}\N_{r,s}.$$

Following Baranovsky, we have that $\overline\N_{0,n}=\N'_{0,n}$
and $\overline\N_{n,0}=\N'_{n,0}$ which can be identified with
$\mathcal{K}_n\times V$ and $\mathcal{K}_n\times V^*$,
respectively.

Our main theorem is the following.
\begin{thm}\label{components}
\begin{enumerate}
\item[(a)] The irreducible components of $\N$ are precisely
$\overline\N_{t,n-1-t},\,1\leq t\leq n-2,\,\N'_{0,n}$ and
$\N'_{n,0}$. \item[(b)] $\dim\overline\N_{t,n-1-t}=n^2+n-2$ for
$1\leq t\leq n-2$ and $\dim \N'_{0,n}=\dim\N'_{n,0}=n^2+n-1$.
\end{enumerate}
\end{thm}

In the last section we study the variety
$$\S=\{(A,B,i,j)\in\n\times\n\times V\times V^*\,|\,A+B=ij\}.$$
and describe its irreducible components. The key fact is that if
$(A,B,i,j)\in \S$ then $A$ and $B$ are simultaneously
triangularizable. The considerations that we make for this variety are simpler that the ones for $\N$ since in this case it is easy to deform an element of $\S$ and stay inside $\S$, one just has to consider matrices that are strictly upper triangular in the given basis.

\section{The Hilbert scheme}\label{Hilbert}

In this section we establish a connection between the Hilbert
scheme and $\N_{t,n-1-t}$ to prove that the later is irreducible
if ${\rm char}\C=0$ or $\geq n/2$. This will also allow us to
prove part (b) of Theorem \ref{components}.

%\begin{defi}
%Let $(X,Y,i,j)\in\N$.  $\C\langle X,Y\rangle i$ (resp.\
%$j\C\langle X,Y\rangle$) is the smallest subspace of $V$ (resp.\
%$V^*$) that contains $i$ (resp.\ $j$) and is invariant under the
%actions of $X$ and $Y$.
%\end{defi}

Clearly every element of $\C\langle X,Y\rangle i$ (resp.\
$j\C\langle X,Y\rangle$) can be written as $p(X,Y)i$ (resp.\
$jp(X,Y)$) where $p(x,y)\in \C\langle x,y\rangle$ is a polynomial
in the noncommutative variables $x$ and $y.$

\begin{lemma}\label{commute}
$j\C\langle X,Y\rangle$ (resp.\ $\C\langle X,Y\rangle i$) is a
right (resp.\ left) $\C[x,y]$-module and its perpendicular
complement $(j\C\langle X,Y\rangle)^\perp=\{v\in V\, |\, uv=0\
\forall u\in j\C\langle X,Y\rangle\}$ (resp.\ $(\C\langle
X,Y\rangle i)^\perp=\{z\in V^*\, |\, zw=0\ \forall w\in \C\langle
X,Y\rangle i\}$) is a left (resp.\ right) $\C[x,y]$-module, where
$x$ and $y$ act as $X$ and $Y$ respectively.
\end{lemma}

\begin{proof}
$j\C\langle X,Y\rangle$ and $(j\C\langle X,Y\rangle)^\perp$ are
right and left $\C\langle x,y\rangle$-modules, respectively. We
have to prove that the two-sided ideal generated by $[x,y]$ acts
as zero on both of them.

Let $q_1(x,y),\,q_2(x,y)\in\C\langle x, y\rangle$. Then for every
$v\in (j\C\langle X,Y\rangle)^\perp$,
$$q_1(X,Y)(XY-YX)q_2(X,Y)v=q_1(X,Y)i(jq_2(X,Y)v)=0$$
since $jq_2(X,Y)\in \C\langle X,Y\rangle.$ This proves that
$(j\C\langle X,Y\rangle)^\perp$ is a $\C[x,y]$-module

To prove that $j\C\langle X,Y\rangle$ itself is a
$\C[x,y]$-module, note that $j\C\langle X,Y\rangle\subseteq
(\langle X,Y\rangle i)^\perp$ since the elements of $j\C\langle
X,Y\rangle$ and $\C\langle X,Y\rangle i$ have the form
$jp(X,Y),\,q(X,Y)i$ for $p(x,y),\,q(x,y)\in\C\langle x, y\rangle$,
and $jp(X,Y)q(X,Y)i=0$ (\cite{Gan-Ginzburg:06}, Lemma 2.1.3.).
This proves the claim.
\end{proof}

From now on, we write $j\C[X,Y]=j\C\langle
X,Y\rangle,\,\C[X,Y]i=\C\langle X,Y\rangle i$.

We now define a map $\N_{r,s}\to\H_{[0]}^r\times\H_{[0]}^s$ where
$\H_{[0]}^m$ is the fiber of the Hilbert-Chow morphism
$\H^m(\A^2)\to S^m(\A^2)$ over the point $m\cdot[0]\in S^m(\A^2)$
and $\H^m(\A^2)$ denotes the Hilbert scheme of points in the
affine plane. The image of $(X,Y,i,j)\in\N_{r,s}$ to $\H_{[0]}^r$
is the ideal $\{p(x,y)\in\C[x,y]\,|\,p(X,Y)i=0\}$. This ideal is
also equal to $\{p(x,y)\in\C[x,y]\,|\,p(X,Y)v=0\ \forall
v\in\C[X,Y]i\}$ since $i$ is a cyclic vector for $X$ and $Y$ on
$\C[X,Y]i.$ Similarly, we define $\N_{r,s}\to\H_{[0]}^s$ as the
ideal $\{p(x,y)\in\C[x,y]\,|\,jp(X,Y)=0\}$.

Since $(\C[X,Y]i)^\perp$ and $(j\C[X,Y])^\perp$ are
$\C[x,y]$-modules, one could try to induce maps
$\N_{r,s}\to\H_{[0]}^{n-s},\H_{[0]}^{n-r}$, but there may not be a
cyclic vector for the actions of $X$ and $Y$ on those spaces.
However, if $r+s=n-1$ the following theorem implies that such map
does exist.

\begin{thm}\label{image}
There is a well-defined regular map
$\N_{t,n-1-t}\to\H_{[0]}^{t+1}(\A^2)$ induced by the actions of
$X$ and $Y$ on $(j\C[X,Y])^\perp$. This map is dominant and the
image of any element of $\N_{t,n-1-t}$ has the form $\langle
y^{t+1}, x-a_1y-\dots -a_{t}y^{t}\rangle$ or $\langle x^{t},
y-a_1x-\dots -a_{t}x^{t}\rangle$ for some $a_1,\dots, a_{t}\in\C$.
\end{thm}

Similar considerations hold for the map
$\N_{t,n-1-t}\to\H_{[0]}^{n-t}(\A^2)$ induced by the actions of
$X$ and $Y$ on $(\C[X,Y]i)^\perp$.  In order to prove Theorem
\ref{components}, we study the image and the fibers of the map
$\Psi:\N_{t,n-1-t}\to \H_{[0]}^t\times \H_{[0]}^{n-t}$. Recall
that $GL_n(\C)$ acts on $\N$ by
$G\cdot(X,Y,i,j)=(GXG^{-1},GYG^{-1},Gi,jG^{-1}).$

\begin{thm}\label{orbits}
The fibers of the map $\Psi:\N_{t,n-1-t}\to \H_{[0]}^t\times
\H_{[0]}^{n-1-t}$ are the $GL_n(\C)$-orbits, and the isotropy of
each element of $\N_{t,n-1-t}$ is one-dimensional.
\end{thm}

%\section{The proofs} \label{proofs}
In order to prove Theorems \ref{image} and \ref{orbits} we need a
technical result.

\begin{lemma}\label{basis}
Let $(X,Y,i,j)\in\N_{r,s}$. There is a basis $\{e_1,\,\dots e_n\}$
of $V$ with dual basis $\{e^*_1,\,\dots e^*_n\}$ of $V^*$ so that
$e_r=i,\,e^*_{n+1-s}=j$, and $X,\,Y$ are upper triangular in this
basis.
\end{lemma}

\begin{proof}
Since $j\C[X,Y]$ annihilates $\C[X,Y]i$ we can decompose
$V=V_1\oplus V_2\oplus V_3$ so that $V_1=\C[X,Y]i$ and
$V_3^*=(j\C[X,Y])^*$. We are to find elements $e_1,\dots,\,e_r\in
V_1,\,e_{r+1},\dots,\,e_{n-s}\in V_2,\,e_{n+1-s},\dots,\,e_n\in
V_3$ satisfying the conditions.

Consider the lex deg order of the monomials in $\C[x,y]$
$$1<x<y<x^2<xy<y^2<x^3<x^2y<xy^2<y^3<\dots$$
Choose the largest (according to $<$) monomial $m_1$ so that
$m_1(X,Y)i\neq0$ (it exists since $X^aY^bi=0$ if $a+b\geq n$), and
inductively choose $m_k$ as the largest monomial so that
$m_k(X,Y)i$ is not a linear combination of
$m_{1}(X,Y)i,\,\cdots,\,m_{k-1}(X,Y)i$. This gives us $r$
monomials $m_1,\,\cdots,\,m_r$ so that $\C[X,Y]i=\langle
m_1(X,Y)i,\,\cdots,\,m_r(X,Y)i\rangle$. We set
$e_1=m_1(X,Y)i,\dots,\,e_r=m_r(X,Y)i=i$. The action of $X$ and $Y$
in this basis is triangular, since multiplying by $x$ or $y$
``increases" monomials, and for every monomial $m$, either
$m(X,Y)i$ is in the basis or is a linear combination of larger
monomials. That $e_r=i$ is a consequence of the following lemma
which will also be useful later.

\begin{lemma}\label{Young} Let $m_t(x,y)=x^{a_t}y^{b_t}$ and
$\lambda=\{(a_1,b_1),\dots,(a_r,b_r)\}$. If $(a,b)\in \lambda$ and
$0\leq a'\leq a,\,0\leq b'\leq b$ then $(a',b')\in \lambda$.
\end{lemma}

\begin{proof}
If $(a',b')\notin \lambda$ then $X^{a'}Y^{b'}i$ can be written as
a linear combination of larger monomials
$\displaystyle\sum_{x^{a'}y^{b'}<x^cy^d}\alpha_{c,d}X^cY^di.$ Then
$$m_t(X,Y)i=X^{a-a'}Y^{b-b'}\(\sum_{x^{a'}y^{b'}<x^cy^d}\alpha_{c,d}X^cY^di\)=
\sum_{x^{a'}y^{b'}<x^cy^d}\alpha_{c,d}X^{c+a-a'}Y^{d+b-b'}i$$ but
this is a contradiction since $x^{a'}y^{b'}<x^cy^d$ implies
$x^{a}y^{b}<x^{c+a-a'}y^{d+b-b'}.$
\end{proof}
In particular this implies that $m_r(x,y)=1$ and therefore
$e_r=i$.

Now we can follow the same procedure in $V_3^*=j\C[X,Y]$ to find
elements $e_{n+1-s}^*,\dots,\,e_n^*\in j\C[X,Y]$ which in turn
give rise to $e_{n+1-s},\dots,\,e_n\in V_3\subseteq V.$

To find the remaining elements of the basis, note that since
$V_2\subseteq(j\C[X,Y])^\perp$ we have that $j\C[X,Y]X
V_2,j\C[X,Y]YV_2\subseteq j\C[X,Y]V_2=0.$ This means that
$X|_{V_2},Y|_{V_2}:V_2\to V_1\oplus V_2$. Since $X,Y:V_1\to V_1$,
the actions of $X$ and $Y$ as endomorphisms of $V_2$ are
nilpotent, and they are commutative since $[X,Y]=ij$ which acts as
0 on this space. Therefore we can find $e_{r+1},\dots,\,e_{n-s}\in
V_2$ that make both $X$ and $Y$ upper triangular when restricted
to $V_2.$ The basis $\{e_1,\,\dots e_n\}$ of $V$ satisfies the
conditions.
\end{proof}

Therefore we can assume that $X$ and $Y$ are upper triangular,
$$i=\begin{pmatrix}
0\\0\\\vdots\\0\\1\\0\\\vdots\\0\\0
\end{pmatrix},\ \C[X,Y]i=\begin{pmatrix}
*\\ *\\\vdots\\ *\\ *\\ 0\\\vdots\\0\\0
\end{pmatrix},j=\begin{pmatrix}
0&0&\dots&0&1&0&\dots&0&0
\end{pmatrix},$$
$$\ j\C[X,Y]=\begin{pmatrix}
0&0&\dots&0&*&*&\dots&*&*
\end{pmatrix};$$
where the 1's in $i$ and $j$ are located in the $r$-th and
$(n+1-s)$-th position, respectively.

\begin{proof}[Proof of Theorem \ref{image}]
We have to prove that $(j\C[X,Y])^\perp$ admits a cyclic vector
and that one of $X,\,Y$ is regular when restricted to that space.

Let $w\in(j\C[X,Y])^\perp\setminus\C[X,Y]i$. Since $Xw,Yw\in
(j\C[X,Y])^\perp=\C w\oplus\C[X,Y]i$ and $X,Y$ are nilpotent, we
have that $Xw,Yw\in\C[X,Y]i$; let $Xw=P(X,Y)i,Yw=Q(X,Y)i,$ where
$P$ and $Q$ are polynomials. Let $Q(x,y)=Q_1(x,y)+c$ where
$Q_1(0,0)=0$. If $c\neq0$ then we have
$$\frac1c (YP(X,Y)-XQ_1(X,Y))i=\frac1c (YP(X,Y)-XQ(X,Y)+cX)i$$
$$=\frac1c(YX-XY)w+Xi=\frac1c ijw+Xi=Xi$$ and every monomial in the left-hand
side is $>x$ in the lex deg order. Multiplying by $X^{a-1}Y^b,\,a>
0,\, b\geq0$; we conclude that $X^aY^bi$ is a linear combination
of larger monomials. According to the construction of the basis in
Lemma \ref{basis} we have that the basis for $\C[X,Y]i$ is
$\{Y^ti,\,\dots,\,Yi,\,i\}$ and this implies that $Y$ acts
regularly in $\C[X,Y]i$.

If $P$ has a constant term we can reverse the roles of $x,y$ in
the lex deg order. Now we prove that at least one of the
polynomials $P,Q$ has a nonzero constant term.

We can choose $u\in V\setminus(j\C[X,Y])^\perp$ so that
$Xu,Yu\in(j\C[X,Y])^\perp$: to do this take any $u_0\in
V\setminus(j\C[X,Y])^\perp$; if $Xu_0,Yu_0\in(j\C[X,Y])^\perp$,
take $u=u_0$. If not, say $Xu_0\notin(j\C[X,Y])^\perp$, take
$u_1=Xu_0$ and repeat the process.

Then $ju \neq 0$ and $Xu=\alpha w+R(X,Y)i, Yu=\beta w+S(X,Y)i$ for
some $R,S\in \C[x,y],\,\alpha,\beta\in\C$. If $ju=0$ then
$j\C[X,Y]u=0$ since $Xu, Yu\in(j\C[X,Y])^\perp$. We can assume,
normalizing $u$ if necessary, that $ju=1$. Therefore
$$i=iju=XYu-YXu=X(\beta w+S(X,Y)i)-Y(\alpha w+R(X,Y)i)$$
$$=(\beta P(X,Y)-\alpha Q(X,Y)+XS(X,Y)-YR(X,Y))i.$$
If $P$ and $Q$ have no constant term then all the monomials in the
right-hand side have positive degree, a contradiction.

Therefore one of $X$ or $Y$ acts regularly on $\C[X,Y].$ Assume
without loss of generality that it is $Y$. It is easy to see that
then $X|_{\C[X,Y]i}=A(Y)|_{\C[X,Y]i}$ where
$A(y)=a_1+\dots+a_{t-1}y^{t-1}\in\C[y]$

Now we prove that there exists
$i'\in(j\C[X,Y])^\perp\setminus\C[X,Y]i$ so that $Yi'=i$. Since
$Xw,Yw\in\C[X,Y]i=\C[Y]i$, let
$$Xw=D(Y)i=(d_0+d_1Y+\dots+d_{t-1}Y^{t-1}i),\,Yw=C(Y)i=(c_0+c_1Y+\dots+c_{t-1}Y^{t-1})i.$$
If $c_0\neq0$ let
$i'=\frac1{c_0}(w-c_1i-\dots-c_{t-1}Y^{t-2}i)\in\C
w\oplus\C[Y]i=(j\C[X,Y])^\perp$. So assume by contradiction that
$c_0=0.$ Then $jw=0$ implies
$$0=ijw=(XY-YX)w=(A(Y)C(Y)-YD(Y))i$$
and therefore, comparing the coefficient of $Y$, $d_0=0$.

Let $u\in V$ as before. Then
$$i=iju=XYu-YXu=X(\beta w+S(X,Y)i)-Y(\alpha w+R(X,Y)i)$$
$$=(\beta D(Y)-\alpha C(Y)+XS(X,Y)-YR(X,Y))i;$$
but all the monomials on the right-hand side have positive degree,
which is impossible. Therefore we can find such $i'.$

But
$Xi'=\frac1{c_0}(Xw-c_1Xi-\dots-c_{t-1}XY^{t-2}i)\in\C[X,Y]i=\subseteq(j\C[X,Y])^\perp.$
Therefore $\C[X,Y]i\subsetneq\C[X,Y]i'\subseteq(j\C[X,Y])^\perp$.
This and $\dim(j\C[X,Y])^\perp=1+\dim\C[X,Y]i$ imply that
$(j\C[X,Y])^\perp=\C[X,Y]i'=\C[Y]i'$ which means that $i'$ is a
cyclic vector in $(j\C[X,Y])^\perp$ and therefore we have a map
$\N_{t,n-1-t}\to\H_{[0]}^{t+1}(\A^2)$. If
$Xi'=(a_1Y+a_2Y^2+\dots+a_{t}Y^{t})i'$ then the image of
$(X,Y,i,j)$ in $\H_{[0]}^{t+1}(\A)$ is $\langle y^{t+1},
x-a_1y-a_2y^2-\dots-a_{t}y^{t}\rangle$.
\end{proof}

\begin{proof}[Proof of Theorem \ref{orbits}] Let $x_{r,s}$ and $y_{r,s}$ denote the entries in the $r$-th row
and $s$-th column of the matrices that represent $X$ and $Y$
respectively in the basis described in Lemma \ref{basis}.

Since
$\begin{vmatrix}x_{t+1,t}&y_{t+1,t}\\x_{t+2,t+1}&y_{t+2,t+1}\end{vmatrix}\neq0$
we can assume without loss of generality that
$x_{t+1,t},y_{t+2,t+1}\neq0$. Therefore
$$\C[X,Y]i=\langle
i,Yi,\dots,Y^{t-1}\rangle,\,j\C[X,Y]=\langle
j,jX,\dots,jX^{n-1-t}\rangle.$$  Consider the filtration
$W_0=(j\C[X,Y])^\perp=\langle
j,jX,\dots,jX^{n-1-t}\rangle^\perp,\,W_1=$\\ $\langle
jX,\dots,jX^{n-1-t}\rangle^\perp,\,\dots,W_{n-1-t}=(jX^{n-1-t})^\perp,$

We choose a basis $\{v_1,\dots,v_t,v_{t+1},\dots,v_n\}$ of $V$ so
that $v_1=Y^{t-1}i,\dots v_{t}=i,v_{t+1}=i'\in
W_0\setminus\C[X,Y]i,\,v_{t+2}\in W_1\setminus W_0,\dots,\ v_n\in
V\setminus W_{n-1-t}$, and we can do this in such a way that
$$x_{p,q}=\begin{cases}
0& q-p\geq 0\\a_1& q-p=1,\,1\leq p\leq t\\1& q-p=1,\,t+1\leq p\leq n-1\\0& q-p\neq1,\,t+1\leq p\leq n,\,t+1\leq q\leq n\\
\end{cases},$$
$$y_{p,q}=\begin{cases}
0& q-p\geq 0\\1& q-p=1,\,1\leq p\leq t\\b_1& q-p=1,\,t+1\leq p\leq n-1\\0& q-p\neq1,\,1\leq p\leq t+1,\,1\leq q\leq t+1\\
\end{cases};$$
where $\Psi(X,Y,i,j)=(y^{t+1},\langle x-a_1y-\dots-a_t y^t\rangle,
\langle x^{n-t},y-b_1x-\dots-b_{n-1-t} x^{n-1-t}\rangle)$ and
$a_1b_1\neq1$

For example, for $n=7,\,t=4$:
$$X=\begin{pmatrix}
0&a_{1}&a_{2}&a_{3}&a_{4}&x_{16}&x_{17}\\0&0&a_{1}&a_{2}&a_{3}&x_{26}&x_{27}\\0&0&0&a_{1}&a_{2}&x_{36}&x_{37}
\\0&0&0&0&a_{1}&x_{46}&x_{47}\\0&0&0&0&0&1&0\\0&0&0&0&0&0&1\\0&0&0&0&0&0&0
\end{pmatrix},\ Y=\begin{pmatrix}
0&1&0&0&0&y_{16}&y_{17}\\0&0&1&0&0&y_{26}&y_{27}\\0&0&0&1&0&y_{36}&y_{37}
\\0&0&0&0&1&y_{46}&y_{47}\\0&0&0&0&0&b_1&b_2\\0&0&0&0&0&0&b_1\\0&0&0&0&0&0&0
\end{pmatrix},$$
$$i=\begin{pmatrix}0\\0\\0\\1\\0\\0\\0\end{pmatrix},i'=\begin{pmatrix}0\\0\\0\\0\\1\\0\\0\end{pmatrix},
j=\begin{pmatrix}0&0&0&0&0&a_1b_1-1&0\end{pmatrix}.$$
%$\Psi(X,Y,i,j)=(\langle x-a_1y-a_2y^2-a_3y^3-a_4 y^4\rangle,
%\langle y-b_1x-b_2 x^2\rangle);\, a_1b_1\neq1$
(Note that we are changing $j$ from our previous notation, this is
to simplify our expressions for $X$ and $Y$).

First we are to prove that the isotropy is one-dimensional. Let
$Z=G-I$ where $G\cdot(X,Y,i,j)=(X,Y,i,j).$ $Z\C[Y]i=0$ and
$j\C[X]Z=0$ imply $z_{p,q}=0$ if $p\geq t+2$ or $q\leq t$. Therefore $Z$ is upper triangular. We want to prove that $z_{p,q}=0$ unless $p=1,\,q=n.$ We proceed by
induction on $q-p.$

$[Z,X]=[Z,Y]=0$ imply
$$\sum_{r\leq m\leq s}\begin{vmatrix}x_{r,m}&z_{r,m}\\x_{m,s}&z_{m,s}\end{vmatrix}=
\sum_{r\leq m\leq
s}\begin{vmatrix}y_{r,m}&z_{r,m}\\y_{m,s}&z_{m,s}\end{vmatrix}=0$$
for every $r\leq s$ (see \cite{Basili2}).

For $r=t, s=t+1$ we have
$$0=\begin{vmatrix}y_{t,t}&z_{t,t}\\y_{t,t+1}&z_{t,t+1}\end{vmatrix}
+\begin{vmatrix}y_{t,t+1}&z_{t,t+1}\\y_{t+1,t+1}&z_{t+1,t+1}\end{vmatrix}=
\begin{vmatrix}0&0\\1&0\end{vmatrix}
+\begin{vmatrix}1&0\\0&z_{t+1,t+1}\end{vmatrix}=z_{t+1,t+1}.$$
This proves the case $q-p=0$.

Assume that $z_{p,q}=0$ if $q-p<u$ and let $r$ and $s$ be
so that $s-r=u+1$. We are to prove that $z_{r,s-1}=z_{r+1,s}=0$.

$$0=\sum_{r\leq m\leq s}\begin{vmatrix}x_{r,m}&z_{r,m}\\x_{m,s}&z_{m,s}\end{vmatrix}=
\begin{vmatrix}x_{r,r}&z_{r,r}\\x_{r,s}&z_{r,s}\end{vmatrix}+
\begin{vmatrix}x_{r,r+1}&z_{r,r+1}\\x_{r+1,s}&z_{r+1,s}\end{vmatrix}
+\begin{vmatrix}x_{r,s-1}&z_{r,s-1}\\x_{s-1,s}&z_{s-1,s}\end{vmatrix}+
\begin{vmatrix}x_{r,s}&z_{r,s}\\x_{s,s}&z_{s,s}\end{vmatrix}$$
$$=\begin{vmatrix}x_{r,r+1}&0\\x_{r+1,s}&z_{r+1,s}\end{vmatrix}
+\begin{vmatrix}x_{r,s-1}&z_{r,s-1}\\x_{s-1,s}&0\end{vmatrix}=
x_{r,r+1}z_{r+1,s}-x_{s-1,s}z_{r,s-1}
=\begin{vmatrix}x_{r,r+1}&z_{r,s-1}\\x_{s-1,s}&z_{r+1,s}\end{vmatrix}.$$
Similarly, $0=\begin{vmatrix}y_{r,r+1}&z_{r,s-1}\\y_{s-1,s}&z_{r+1,s}\end{vmatrix}$.

$r\leq t$ and $s-1>t$ we have
$$0=\begin{vmatrix}a_1&z_{r,s-1}\\1&z_{r+1,s}\end{vmatrix}=
\begin{vmatrix}1&z_{r,s-1}\\b_1&z_{r+1,s}\end{vmatrix}\Rightarrow z_{r,s-1}=z_{r+1,s}=0.$$

If $r>t$ then $r+1\geq t+2$, so $z_{r+1,s}=0$ and
$$0=\begin{vmatrix}x_{r,r+1}&z_{r,s-1}\\x_{s-1,s}&z_{r+1,s}\end{vmatrix}=
\begin{vmatrix}1&z_{r,s-1}\\1&0\end{vmatrix}=-z_{r,s-1}.$$
Similar considerations apply if $s-1\leq t.$

This proves that $z_{p,q}=0$ if $p-q\leq n-2$ (since in the
computations above we require $r\geq 1,\,s\leq n$). Therefore
$z_{p,q}=0$ unless $p=1,\,q=n$.

It is easy to check that if $z_{1,n}$ is arbitrary and $z_{p,q}=0$
for $(p,q)\neq(1,n)$ then $G=I+Z$ fixes $(X,Y,i,j).$ Therefore the
isotropy is one-dimensional.

Now we want to prove that the conjugacy class is uniquely determined
by the image under $\Psi.$ In order to do that we are to find a suitable basis in which $X$ and $Y$ are easy to describe.

\begin{lemma}\label{vector}
There exists a vector $v\in V$ so that
\begin{enumerate}
\item[(a)] $v\notin Im\,X + Im\,Y$ ; \item[(b)] $Z_mv\in Y\C[X,
Y]i$ where $Z_m= (Y-b_1X-\dots-b_{n-1-t-m}X^{n-1-t-m})X^m; 0\leq m
\leq n-2-t.$ \item[(c)] $Z_0v = (Y-b1X-\dots-b_{n-1-t}X^{n-1-t})v
= 0$;
\end{enumerate}
\end{lemma}

\begin{proof} We will construct $v$ one entry at a time (using the basis from before).
$$v=\begin{pmatrix}v_1\\ \vdots\\ v_{t+1}\\ v_{t+2}\\ \vdots\\v_n\end{pmatrix}.$$

Let $v_n=1$, this guarantees (a).

Now we prove that $Im\,Z_m \subseteq\C[X,Y]i$. This means that the
last $n-t$ rows of $Z_m$ vanish. Let $\tilde Y$ and $\tilde X$ be
the $(n-t)\times(n -t)$ lower right submatrices of $Y$ and $X$
respectively; so in fact $\tilde X$ is a regular nilpotent matrix,
and $\tilde Y=b_1\tilde X+\dots + b_{n-1-t}\tilde X^{n-1-t}$.
Therefore
$$(\tilde Y-b_1 \tilde X -\dots- b_{n-1-t-m}\tilde X^{n-1-t-m})\tilde X^m =$$
$$(b_{n-t-m}\tilde X^{n-t-m} + \dots+ b_{n-1-t}\tilde X^{n-1-t})\tilde X^m = 0.$$

It is easy to see that the $(t+1)\times(t+1)$ upper left submatrix
of $Y-b_1X-\dots -b_{n-1-t-m}X^{n-1-t-m}$ is a regular nilpotent
matrix, in fact, the entries just above the main diagonal are all
equal to $1-a_1b_1\neq0$.

Now consider the $t$-th row of $Z_m$. The $(t+1)$-th entry of that
row of $(Y- b_1X-\dots- b_{n-1-t-m}X^{n-1-t-m}$ is $1-a_1b_1$, and
all the preceding entries are equal to 0. After multiplying by
$X^{n-1-t-m}$, the $(t+1+m)$-th entry of the $t$-th row of $Z_m$
is equal to $1-a_1b_1$, and all the preceding entries are equal to
0. So we use $Z_{n-2-t},\,Z_{n-3-t},\dots,Z_1$ to choose
$v_{n-1},v_{n-2},\dots,v_{t+2}$ so that the $t$-th entry of $Z_mv$
is zero, i.e., $Z_mv\in Y \C[X,Y]i =Y\C[Y]i =\langle Yi,
Y^2i,\dots,Y^ti\rangle$. 

Since the $(t+1)\times(t+1)$ upper left submatrix of $Z_0$ is
regular then, given $v_{t+2},\dots,v_n\in\C$, there are unique
numbers $v_2,\dots,v_{t+1}\in\C$ so that the vector $w$ formed in
this way is in $\ker Z_0$ ($w_1$ is arbitrary).

\end{proof}

Since $v\in V\setminus Im\,X+Im\,Y$ we have that $v$ is a cyclic
vector for $X$ and $Y$, in fact,
$$v=\begin{pmatrix}*\\ *\\\vdots\\ *\\ *\\1\end{pmatrix},\,
Xv=\begin{pmatrix}*\\ *\\\vdots\\ *\\ 1\\0\end{pmatrix},\,
X^2v=\begin{pmatrix}*\\ *\\\vdots\\ 1\\ 0\\0\end{pmatrix},\dots,\,
X^{n-1-t}v=\begin{pmatrix}*\\\vdots\\ *\\1\\
\vdots\\0\end{pmatrix},$$
$$YX^{n-1-t}v=\begin{pmatrix}*\\\vdots\\1\\0\\
\vdots\\0\end{pmatrix},\dots\,,Y^tX^{n-1-t}v=\begin{pmatrix}1\\\vdots\\
0\\0\\ \vdots\\0\end{pmatrix}.
$$

Therefore $V=\langle v,Xv,\dots,X^{n-1-t}v,YX^{n-1-t}v,\dots,
Y^tX^{n-1-t}v\rangle,\,\C[X,Y]i=\langle YX^{n-1-t}v,\dots,
Y^tX^{n-1-t}v\rangle,\,(j\C[X,Y])^\perp=\langle
X^{n-1-t}v,YX^{n-1-t}v,\dots, Y^tX^{n-1-t}v\rangle$

Now we write $X$ and $Y$ as matrices in this basis. Since $Y$ acts
regularly in $(j\C[X,Y])^\perp$ and $X$ acts as
$a_1Y+\dots+a_{t+1}Y^{t+1}$ in this space we have that
$$x_{p,q}=\begin{cases}
%0& q-p\geq 0\\
a_{q-p}& 1\leq p<q\leq t\\1& q-p=1,\,t+1\leq p\leq n-1\\
0& \text{otherwise}\\
\end{cases}.$$

$v\in\ker Z_0$ means
$$Yv =b_1Xv -\dots- b_{n-1-t}X^{n-1-t}v$$
and for $m = 0,\dots, n-2-t$; the coefficient of $Y X^{n-1-t}v$ in
$Y(X^mv)$ is 0. This implies that

$$y_{p,q}=\begin{cases}
1& q-p=1,\,1\leq p\leq t\\
b_{q-p}&t+1\leq p<q\leq n\\
%0& q-p\neq1,\,1\leq p\leq t+1,\,1\leq q\leq t+1\\
0& q-p\geq 0,\,p=t\text{ or }q=n\\
\end{cases}.$$

For example, for $n=8,\,t=4$:
$$X=\begin{pmatrix}
0&a_{1}&a_{2}&a_{3}&a_{4}&0&0&0\\0&0&a_{1}&a_{2}&a_{3}&0&0&0\\0&0&0&a_{1}&a_{2}&0&0&0
\\0&0&0&0&a_{1}&0&0&0\\0&0&0&0&0&1&0&0\\0&0&0&0&0&0&1&0\\0&0&0&0&0&0&0&1
\\0&0&0&0&0&0&0&0
\end{pmatrix},\ Y=\begin{pmatrix}
0&1&0&0&0&y_{16}&y_{17}&0\\0&0&1&0&0&y_{26}&y_{27}&0\\0&0&0&1&0&y_{36}&y_{37}&0
\\0&0&0&0&1&0&0&0\\0&0&0&0&0&b_1&b_2&b_3\\0&0&0&0&0&0&b_1&b_2\\0&0&0&0&0&0&0&b_1
\\0&0&0&0&0&0&0&0
\end{pmatrix},$$

To conclude the proof of the Theorem \ref{orbits} we use following
\begin{lemma}
The entries $y_{p,q}$ for $1\leq p<t,\,t+2\leq q<n$, the vector $i$ and
the covector $j$ are uniquely determined by $a_1,\dots,\,
a_{t},\,b_1,\dots,\,b_{n-1-t}$ and by the condition $rk([X,Y])=1.$
\end{lemma}

\begin{proof}
We split $X$ and $Y$ in blocks of sizes $t,\,1,\,n-1-t$:
$$X=\begin{pmatrix}X_{11}&X_{12}&0\\0&0&X_{23}\\0&0&X_{33}\\
\end{pmatrix},\,Y=\begin{pmatrix}Y_{11}&Y_{12}&Y_{13}\\0&0&Y_{23}\\0&0&Y_{33}\\
\end{pmatrix},$$
so
$$X_{11}=\begin{pmatrix}0&a_1&\dots&a_{t-2}&a_{t-1}\\
0&0&\dots&a_{t-3}&a_{t-2}\\
\vdots&&&&\vdots\\0&0&\dots&0&a_1
\\0&0&\dots&0&0\end{pmatrix},\,X_{12}=\begin{pmatrix}a_{t}\\a_{t-1}\\
\vdots\\a_2\\a_1 \end{pmatrix}\text{, etc.}$$

Then $X_{11}Y_{13}+X_{12}Y_{23}-Y_{12}X_{23}-Y_{13}X_{33}=\tilde
i\tilde j$ ,a rank one matrix, where $\tilde i$ and $\tilde j$ are
the truncated vector and covector respectively. Note that the
entry in the lower left position of $\tilde i\tilde j$ is equal to
$a_1b_1-1\neq0$, and that $X_{12}Y_{23}$ and $Y_{12}X_{23}$ are
determined by $a_1,\dots,\, a_{t},\,b_1,\dots,\,b_{n-1-t}$.

Let $$h_1=\begin{pmatrix}1&0&\dots&0&0\end{pmatrix},\,
h_2=\begin{pmatrix}0&1&\dots&0&0\end{pmatrix},\dots,\,
h_t=\begin{pmatrix}0&0&\dots&0&1\end{pmatrix}.$$ Since
$h_tX_{11}=0$ and $h_tY_{13}=0$ we conclude that $h_t\tilde
i\tilde j=h_t\(X_{12}Y_{23}-Y_{12}X_{23}\)\neq 0$. From here we
get that $\tilde j$ is uniquely determined by the parameters
$a_1,\dots,\, a_{t}$ $b_1,\dots,\,b_{n-1-t}$ (up to a constant
multiple, that we ignore).

In general, after we know $h_mY_{13}$ for $m> k$ , we calculate
$$h_{k}\tilde i\tilde
j=h_{k}\(X_{11}Y_{13}+X_{12}Y_{23}-Y_{12}X_{23}-Y_{13}X_{33}\).$$
Since $h_kX_{11}=a_1h_{k-1}+a_2h_{k-2}+\dots$ and
$h_kY_{13}X_{33}=\begin{pmatrix}0&y_{k,t+2}&\dots&y_{k,n-1}\end{pmatrix}$
we conclude that the first entry of
$h_{k}\(X_{11}Y_{13}+X_{12}Y_{23}\right.$
$\left.-Y_{12}X_{23}-Y_{13}X_{33}\)$ does not depend on $Y_{13}$, so
neither does $h_{k}\tilde i\tilde j$. This means that $h_{k}\tilde
i$ depends only on $a_1,\dots,\, a_{t},\,b_1,\dots,\,b_{n-1-t}$ (the
key fact here is that $a_1b_1-1\neq0$). And with this, we can find
$h_kY_{13}X_{33}=\begin{pmatrix}0&y_{k,t+2}&\dots&y_{k,n-1}\end{pmatrix}$.
Since the last entry of
$h_kY_{13}=\begin{pmatrix}y_{k,t+2}&\dots&y_{k,n-1}&0\end{pmatrix}$
is zero, we are not loosing any information.
\end{proof}
This concludes the proof of Theorem \ref{orbits}.
\end{proof}

From Theorem \ref{orbits} and the fact the $\Psi$ is dominant we
conclude that $\N_{t,n-1-t}$ (and therefore
$\overline\N_{t,n-1-t}$) is irreducible and its dimension is
$\dim\H_{[0]}^{t+1}+\dim\H_{[0]}^{n-t}+\dim GL(V)-1=n^2+n-2$.
$\dim\N_{0,n}=\dim\N'_{n,0}=n^2+n-1$ as was proved in
\cite{Baranovsky:01}.

\section{Proof of Theorem \ref{components}}\label{cons}

First we prove the non-redundancy of $\overline\N_{t,n-1-t}; 1\leq
t\leq n-2.$
\begin{lemma}
$\N_{t,n-1-t}\cap\overline\N_{r,s}=\emptyset$ unless $r=t$ and
$s=n-1-t.$
\end{lemma}
\begin{proof}
Since $\overline\N_{r,s}\subseteq\N'_{r,s}$, then
$\N_{t,n-1-t}\cap\overline\N_{r,s}\neq\emptyset$ implies $t\leq
r,\,n-1-t\leq s$, so $n-1\leq r+s$, but $\N_{r,s}=\emptyset$ if
$r+s\geq n$ and $r,s>0$, therefore $n-1= r+s$ and we conclude that
$t=r,\,n-1-t=s.$
\end{proof}

To complete the proof of Theorem \ref{components} we have to prove
that
$$\N=\N'_{0,n}\cup\N'_{n,0}\cup\bigcup_{t=1}^{n-2}\overline\N_{t,n-1-t}.$$

\begin{thm}\label{closure}
If $0<r+s<n-1$ then $\N_{r,s}\subseteq
\bigcup_{0<t<n-1}\overline\N_{t,n-1-t}.$
%$$\N=\N'_{0,n}\cup\N'_{n,0}\cup\bigcup_{0<t<n-1}\overline\N_{t,n-1-t},$$
\end{thm}

Let $(X,Y,i,j)\in\N,\,r=\dim\C[X,Y]i$. Let $\lambda$ be as in
Lemma \ref{Young}. The conclusion of the Lemma means that
$\lambda$ represents a Young diagram of size $r$. Let
$\lambda_x=\{(a,b)\in \lambda|(a+1,b)\in
\lambda\},\,\lambda_y=\{(a,b)\in \lambda|(a,b+1)\in \lambda\}$.

Fix a decomposition $V=\C[X,Y]i\oplus V'$ and write $X$ and $Y$ in blocks accordingly:
$$X=\begin{pmatrix} X_1&X_3\\0&X_2 \end{pmatrix},\,Y=\begin{pmatrix}
Y_1&Y_3\\0&Y_2\end{pmatrix}$$
%,\,i=\begin{pmatrix},\, i_1\\0\end{pmatrix}.$$
Here $X_1$ and $Y_1$ (resp.\ $X_3$ and $Y_3$) are nilpotent
commuting endomorphisms of $\C[X,Y]i$ (resp.\ $V'$) while
$X_3,Y_3:V'\to\C[X,Y]i.$ We want to describe the subvariety
\begin{multline*}
\N_{X_1,Y_1,i}=\left\{(X',Y',j')\in \n\times\n\times V^*|\right. \\\left. X'
=\begin{pmatrix} X_1&X'_3\\0&X'_2
\end{pmatrix},\,Y'=\begin{pmatrix} Y_1&Y'_3\\0&Y'_2\end{pmatrix},[X',Y']=ij'\right\}
\end{multline*}

\begin{lemma}\label{birational1}
$\N_{X_1,Y_1,i}$ is birationally equivalent to $\K_{n-r}\times
(V'^*)^{r+1}$ and therefore irreducible.
\end{lemma}
In concrete terms, we have the following
\begin{lemma}\label{birational2}
Let $$X'_3=\sum_{(a,b)\in
\lambda}X^{a}Y^{b}i\alpha_{(a,b)},\,Y'_3=\sum_{(a,b)\in
\lambda}X^{a}Y^{b}i\beta_{(a,b)}$$ where
$\alpha_{(a,b)},\,\beta_{(a,b)}\in V'^*.$ Then
$j',\,\{\beta_{(a,b)}\,|\,(a,b)\in \lambda_x\}$ and
$\{\alpha_{(0,b)}\,|\, (0,b)\in \lambda_y\}$ are uniquely
determined by $X_2',\,Y_2'$ and the remaining
$\alpha_{(a,b)},\,\beta_{(a,b)}.$
\end{lemma}

\begin{proof}
The condition $ij'=X_1Y'_3+Y'_3X'_2-Y_1X'_3-X'_3Y'_2$ means
\begin{multline}\label{eq1}
ij'=\sum_{(a,b)\in
\lambda}(X^{a+1}Y^{b}i)(\beta_{(a,b)})+(X^{a}Y^{b}i)(\beta_{(a,b)}X'_2)-
(X^{a}Y^{b+1}i)(\alpha_{(a,b)})\\-(X^{a}Y^{b}i)(\alpha_{(a,b)}Y'_2)
\end{multline}

\begin{multline*}
=\sum_{(a,b)\in
\lambda}(X^{a}Y^{b}i)(\beta_{(a,b)}X'_2-\alpha_{(a,b)}Y'_2)+
\sum_{(a,b)\in
\lambda,a>0}(X^{a}Y^{b}i)(\beta_{(a-1,b)})-\\\sum_{(a,b)\in
\lambda,b>0}(X^{a}Y^{b}i)(\alpha_{(a,b-1)})+W\end{multline*}
where $W=\sum_{(a,b)\in \lambda\setminus
\lambda_x}(X^{a+1}Y^{b}i)(\beta_{(a,b)})-\sum_{(a,b)\in
\lambda\setminus \lambda_y}(X^{a}Y^{b+1}i)(\alpha_{(a,b)})$. Note
that $W$ depends on $i,\,X_1,\,Y_1,\,\{\beta_{(a,b)}\,|\,(a,b)\in
\lambda\setminus \lambda_x\}$ and $\{\alpha_{(a,b)}\,|\,(a,b)\in
\lambda\setminus \lambda_y\}$. Let $W=\sum_{(a,b)\in \lambda}
(X^{a}Y^{b}i)(w_{(a,b)})$.

Let $(a,b)\in \lambda$. Consider the appearances of $X^{a}Y^{b}i$
in equation \eqref{eq1}. For $(a,b)=(0,0)$, we have
$j'=\beta_{(0,0)}X_2'-\alpha_{(0,0)}Y_2'$, and for every other
$(a,b)\in \lambda,$
\begin{equation}
0=\beta_{(a,b)}X'_2-\alpha_{(a,b)}Y'_2+\beta_{(a-1,b)}-\alpha_{(a,b-1)}+w_{(a,b)}\label{eq2}\end{equation}
where $\beta_{(-1,b)}=\alpha_{(a,-1)}=0$.

Therefore we can define
$$\beta_{(a,b)}=-\beta_{(a+1,b)}X'_2+\alpha_{(a+1,b)}Y'_2+\alpha_{(a+1,b-1)}-w_{(a+1,b)},\quad(a,b)\in \lambda_x;$$
$$\alpha_{(0,b)}=\beta_{(0,b+1)}X'_2-\alpha_{(0,b+1)}Y'_2+w_{(0,b+1)},\quad(0,b)\in \lambda_x.$$
This guarantees \eqref{eq2} for every $(a,b)\in \lambda\setminus
(0,0).$ The definition is non-recursive since every
$\alpha_{(a,b)},\,\beta_{(a,b)}$ is defined as a regular function
of covectors associated to higher monomials in the lex deg order.
\end{proof}

A particular case of Lemma \ref{birational2} is important in its
own right.

\begin{lemma}\label{regularity}
Let $X_1\in\mathfrak{gl}_r,\,Y_2\in\mathfrak{gl}_{n-r}$ be
nilpotent and regular. Then for generic $\alpha_{(a,b)},\,\beta_{(a,b)}$ we have
that $\dim j'\C[X,Y]=n-1-r.$
\end{lemma}

\begin{proof}
Since $X_1$ is regular then $Y_1=\sum_{u=1}^{r-1}c_uX_1^u$ with
$c_1,\dots,c_{r-1}\in\C$, $i$ is a cyclic vector for $X_1$ and
$\lambda=\{(0,0),\dots,({r-1,0})\}.$
%$\lambda=(1,1,\dots),\,\mu=(n)$.
For simplicity we denote
$\alpha_t=\alpha_{(t,0)},\,\beta_t=\beta_{(t,0)}$ so
$$X'_3=\sum_{t=0}^{r-1}X^ti\alpha_{t},\,Y'_3=\sum_{t=0}^{r-1}X^ti\beta_{t}.$$
It follows that
$$ij'=X_1Y'_3+Y'_3X'_2-Y_1X'_3-X'_3Y'_2$$
$$=X_1\(\sum_{t=0}^{r-1}X^ti\beta_{t}\)+\(\sum_{t=0}^{r-1}X^ti\beta_{t}\)X'_2-
\(\sum_{u=1}^{r-1}c_uX_1^u\)\(\sum_{t=0}^{r-1}X^ti\alpha_{t}\)-\(\sum_{t=0}^{r-1}X^ti\alpha_{t}\)Y'_2$$
$$=\sum_{t=0}^{r-1}X^{t}i\(\beta_{t-1}+\beta_{t}X'_2-\alpha_{t}Y'_2-\sum_{u=1}^{t}c_u\alpha_{t-u}\);$$
where we define $\beta_{-1}=0$. Then
$j'=\beta_{0}X'_2-\alpha_{0}Y'_2$ and
$\beta_{t-1}+\beta_{t}X'_2-\alpha_{t}Y'_2-\sum_{u=1}^{t}c_u\alpha_{t-u}=0$
for every $t\geq1.$ This implies that
$\alpha_0,\dots,\alpha_{r-1},\beta_{r-1}\in V'^*$ are arbitrary
and $\beta_{r-2},\dots,\beta_0$ are defined recursively by
$\beta_{t-1}=-\beta_{t}X'_2+\alpha_{t}Y'_2+\sum_{u=1}^{t}c_u\alpha_{t-u}=0$.

Now we want to calculate $\dim j'\C[X',Y']=\dim j'\C[X_2',Y_2']$
in the case where $Y_2'$ is regular and
$X_2'=\sum_{v=1}^{n-r-1}d_vY_2'^v$. Since
$\beta_0=\alpha_0Y_2'+c_1\alpha_0-\beta_1X_2'$, it turns out that
$j'=(c_1d_1-1)\alpha_{0}Y'_2+$terms with higher powers of $Y_2'.$
It follows that $\dim j'\C[X',Y']=n-r-1$ when $c_1d_1-1\neq0$ and
$\alpha_0$ is a cyclic vector for $Y_2'$.
\end{proof}

\begin{proof}[Proof of Theorem \ref{closure}]
Let $(X,Y,i,j)\in\N_{r,s}.$ We want to prove that every open set
$U$ in $\N$ containing $(X,Y,i,j)$ intersects $\N_{t,n-1-t}$ for
some $t$. Since $\N_{X_1,Y_1,i}$ is birationally equivalent to
$\K_{n-r}\times (V'^*)^{r+1}$ and the pairs of nilpotent matrices
where one is regular is dense in $\K_{n-r}$ (see \cite{Basili3}),
we can find $(X',Y',i,j')\in U\cap\N_{X_1,Y_1,i}$ so that $X_2$ is
regular. Now we can reverse the roles of $i$ and $j$ to apply
Lemma \ref{regularity} and find $(X'',Y'',i',j')\in
U\subseteq\N_{t,n-1-t}$ where $t=\dim j'\C[X',Y']$.
\end{proof}

\begin{remark}
The hypothesis about the characteristic of $\C$ has only been used
for the irreducibility of the zero fiber of the Hilbert scheme. In
other characteristics it is still true that
$$\N=\N'_{0,n}\cup\N'_{n,0}\cup\bigcup_{t=1}^{n-2}\overline\N_{t,n-1-t}$$
but the closed sets on the right-hand side may not be irreducible.
\end{remark}

\section{${\rm rk}(A+B)\leq1$}

Consider the variety
$$\S=\{(A,B,i,j)\in\n\times\n\times V\times V^*\,|\,A+B=ij\}.$$

\begin{lemma}\label{triangular}
If $(A,B,i,j)\in S$ then $A$ and $B$ can be simultaneously
triangularized.
\end{lemma}

\begin{proof}
We proceed by induction on $\dim V$, the first case being obvious.
Let $v\in\ker B\setminus\{0\}$. If $jv=0$ then $Av=0$ and we can
apply induction on $V/\C v$. Otherwise we can assume that $jv=1$.
Then $i=ijv=Av+Bv=Av$ and therefore $B=ij-A=A(vj-I)$. Choose $w\in
V^*$ so that $wA=0$. Then $wB=0$ and $wi=0$, so we can apply the
induction on $V^*/w.$
\end{proof}

Recall that $\C\langle A,B\rangle i$ be the smallest subspace of $V$ containing $i$ and invariant under $A$ and $B$.

\begin{lemma}\label{commute1}
$(A+B)|_{\C\langle A,B\rangle i}=0.$
\end{lemma}
The proof of this Lemma is essentially the proof of Lemma 2.1.3 in
\cite{Gan-Ginzburg:06} using our Lemma \ref{triangular}.

\begin{proof}
We have to prove that $(A+B)p(A,B)i=ijp(A,B)i=0$ for any
polynomial $p$ in two noncommutative variables. But
$$jp(A,B)i={\rm tr}(p(A,B)ij)={\rm tr}(p(A,B)(A+B))=0$$
since $A$ and $B$ are upper triangular in the same basis.
\end{proof}

As a consequence we can write $\C[A]i=\C[A,B]i=\C\langle A,B\rangle i$ and similarly $j\C[A]=j\C[A,B]=j\C\langle A,B\rangle.$

Now we want to find the irreducible components of $\S.$ Let
$$\S_{r,s}=\{(A,B,i,j)\in\S\,|\,\dim \C[A]i\leq r,\,\dim j\C[A]\leq s\}.$$
This is clearly a closed set in the Zariski topology.
%and $$\S=\bigcup_{r+s\leq n-1}\S_{r,s}.$$

\begin{thm}
The irreducible components of $\S$ are
$\S_{r,n-r},\,r=0,\,1,\dots,\,n-1.$
\end{thm}

\begin{proof}

Lemma \ref{triangular} implies $\dim \C[A] i+\dim j\C[A]\leq n-1$
since, in the basis given by the Lemma, the matrix $ij$ is upper
triangular and therefore the number of nonzero entries of $i$ and
$j$ cannot exceed $n.$ Therefore
$$\S=\bigcup_{r=0}^{n}\S_{r,n-r}.$$

Let $$L=\begin{pmatrix}0&1&0&\dots&0&0&0\\0&0&1&\dots&0&0&0\\
0&0&0&\dots&0&0&0\\
\vdots&&&\ddots&&\vdots\\0&0&0&\dots&0&1&0
\\0&0&0&\dots&0&0&1\\0&0&0&\dots&0&0&0\end{pmatrix},\,i_r=\begin{pmatrix}0\\
\vdots\\0\\1\\0\\\vdots\\0\end{pmatrix},\,j_r=\begin{pmatrix}0&
\dots&0&1&0&\dots&0\end{pmatrix}$$ where the 1's in $i_r$ and
$j_r$ are located in the $r$-th and $(r+1)$-th positions,
respectively. Then $(L,i_rj_r-L,i_r,j_r)\in\S_{r',n-r'}$ if and
only if $r=r'$. This proves that the closed sets
$\S_{r,n-r},\,r=0,\,1,\dots,\,n-1$ are non-redundant.

Now we have to prove that $\S_{r,n-r}$ is irreducible. Let
$$\S_{r,s}'=\{(A,B,i,j)\in\S\,|\,\dim \C[A]i=r,\,\dim j\C[A]=s\}.$$
We are to prove that $\S_{r,n-r}'$ is irreducible and that its
closure in the Zariski topology is $\S_{r,n-r}.$

Let $(A,B,i,j)\in\S_{r,n-1-r}$ and define $L,\,i_r,\,j_r$ as
before in the same basis that is used to triangularize $A$ and
$B$. For $\tau\in\C$ consider $(A(\tau),B(\tau),i(\tau),j(\tau))$
where $A(\tau)=\tau A+(1-\tau)L,\,i(\tau)=\tau
i+(1-\tau)i_r,\,j(\tau)=\tau
j+(1-\tau)j_r),\,B(\tau)=i(\tau)j(\tau)-A(\tau)$. Clearly the
curve $\{(A(\tau),B(\tau),i(\tau),j(\tau))\,|\,\tau\in\C\}$ is
contained in $\S_{r,n-r}$ and $\dim \C[A(\tau)]i(\tau)=r,\,\dim
j(\tau)\C[A(\tau)]=n-r$ except for finitely many values of $\tau$.
Therefore $(A,B,i,j)$ is in the closure of $\S_{r,n-r}'$.

Every $(A,B,i,j)\in\S_{r,n-r}'$ can be conjugated into
$$A=\begin{pmatrix}J&A'\\0&J
\end{pmatrix},\,i=i_r,\,j=j_r,\,B=ij-A$$ where the blocks in $A$ have
sizes $r$ and $n-r$, respectively, and $J$ represents the Jordan
block of the appropriate size. The elements of the isotropy have
the form $$G=\begin{pmatrix}I&G'\\0&I
\end{pmatrix},\,G'J=JG'.$$
The space of $r\times(n-r)$ blocks $G'$ with $G'J=JG'$ is isomorphic to $\C^{\min\{r,n-r\}}$ and therefore irreducible. We conclude that $\S_{r,n-r}'$ is irreducible and
its dimension is $r(n-r)+n^2-\min\{r,n-r\}.$

\end{proof}

\end{document}